\input amstex
\documentstyle{amsppt}

\magnification=1200
\NoBlackBoxes

\def\th{\theta}
\def\a{\alpha}
\def\l{\lambda}

\def\C{\Bbb C}

\def\Dn{\frak{D}(n;\th)}
\def\Ln{\Lambda^\th(n)}
\def\Pm{P^*_\mu}
\def\O1{{\Cal O}_1}
\def\m1{{\frak m}_1}
\def\S1{{\Cal S}_1}
\def\td{\operatorname{\text{$\th$-}dim}}
\def\G{\Gamma}
\def\ft{)_{\th-1}}
\def\Ts{{\sssize T}(s)}

\def\tht{\thetag}

\topmatter
\title
Shifted Jack polynomials, binomial formula, and applications
\endtitle
\author
A.~Okounkov \footnotemark
and G.~Olshanski \footnotemark
\endauthor
\abstract
In this note we prove an explicit binomial formula for Jack polynomials
and discuss some applications of it.
\endabstract
\thanks
The authors were supported by the Russian Basic Research
Foundation (grant 95-01-00814). The first author was supported 
by the NSF grant DMS 9304580.
\endthanks
\endtopmatter

\leftheadtext{\smc A.~Okounkov and G.~Olshanski}
\rightheadtext{
Shifted Jack polynomials, binomial formula, and applications}

\adjustfootnotemark{-1}
\footnotetext{Institute for Advanced Study, Princeton, NJ 08540.
Current address: Dept.\ of Mathematics, University of
Chicago, 5734 South University Ave., Chicago, IL 60637-1546.
E-mail address: okounkov\@math.uchicago.edu }

\adjustfootnotemark{1}
\footnotetext
{Institute for Problems of Information Transmission,
Bolshoy Karetny 19, 101447 Moscow GSP--4, Russia.
E-mail address: olsh\@ippi.ac.msk.su }

\subhead
1.~Jack polynomials (\cite{M,St})
\endsubhead
In this note we use the parameter
$$
\th=1/\a
$$
inverse to the standard parameter $\a$ for Jack polynomials.
Jack symmetric polynomials $P_\l(x_1,\dots,x_n;\th)$
are eigenfunctions of Sekiguchi differential
operators
$$
\align
&D(u;\th)=V(x)^{-1} \det
\left[
x_i^{n-j}
\left(x_i\frac{\partial}{\partial x_i}+(n-j)\th + u
\right)
\right]_{1\le i,j \le n}\,, \\
&D(u;\th)\, P_\l(x;\th)= 
\left(\prod_i (\l_i+(n-i)\th + u)
\right)\, P_\l(x;\th)\,, \tag1.1
\endalign
$$
where $V(x)=\prod_{i<j}(x_i-x_j)$ is the Vandermonde determinant and
$u$ is a parameter. The operators $\{D(u;\th),u\in\C\}$
generate a commutative algebra; denote it by $\Dn$. 

We normalize $P_\l(x_1,\dots,x_n;\th)$ so that
$$
P_\l(x_1,\dots,x_n;\th)=x_1^{\l_1}\dots x_n^{\l_n} + \dots\,,
$$
where dots stand for lower monomials in lexicographic order.
Then one has 
$$
P_{(\l_1+1,\dots,\l_n+1)} (x;\th)= 
\left(
\prod x_i 
\right)\,P_\l(x;\th)\,. \tag1.2
$$
Using \tht{1.2} one can define Jack rational functions for
any integers $\l_1\ge\dots\ge\l_n$ and still \tht{1.1},
and hence the binomial theorem below, will hold.

Note that the eigenvalues in \tht{1.1} are symmetric
in variables $\l_i-\th i$, $i=1,\dots,n$.

\subhead
2.~Shifted Jack polynomials (\cite{S1,OO,Ok1,KS,Ok3})
\endsubhead
Denote by $\Ln$ the algebra of polynomials $f(x_1,\dots,x_n)$
that are symmetric in variables $x_i-\th i$.
We call such polynomials {\it shifted symmetric}. By \tht{1.1} one
has the Harish-Chandra isomorphism
$$
\Dn\to\Ln \tag2.1
$$
which takes an operator $D\in\Dn$ to the polynomial $d\in\Ln$
such that
$$
D \, P_\l(x;\th)= 
d(\l)\, P_\l(x;\th)\,. 
$$
If $D$ is of order $k$ then $\deg d=k$.

Let $\mu$ be a partition (also viewed as a diagram). 
Recall that for a square  $s=(i,j)\in\mu$ the numbers
$$
\alignat2
&a(s)=\mu_i-j,&\qquad &a'(s)=j-1,\\
&l(s)=\mu'_j-i,&\qquad &l'(s)=i-1,
\endalignat
$$
are called arm-length, arm-colength, leg-length, and
leg-colength respectively. Set
$$
H(\mu)=\prod_{s\in\mu}(a(s)+\th\, l(s)+1) \,.
$$

According to a general result of  S.~Sahi \cite{S1} there exists the
unique polynomial $\Pm(x;\th)$ in 
$\Ln$ such that $\deg \Pm \le|\mu|$ and
$$
\Pm(\lambda;\th)=
\cases H(\mu), &\lambda=\mu\,,\\
0, &|\mu|\le|\lambda|, \; \mu\ne\lambda\,
\endcases
$$
Here we assume $\mu$ and $\lambda$ have length $\le n$.
It is clear that the uniqueness of $\Pm(x,\th)$ follows
easily from its existence.

The polynomials $\Pm(x;1)$ were studied in \cite{OO} and \cite{Ok1}.
They are closely related to Schur
polynomials and are called {\it shifted Schur functions}.
In particular, the Schur function $s_\mu$ is the highest degree term of $\Pm(x;1)$.
Shifted Schur functions have a very rich combinatorics and 
numerous applications (such as Capelli-type 
identities or asymptotic character theory, see \cite{OO,Ok1,N,Ok2} and references
therein). 

For general $\th$, F.~Knop and S.~Sahi proved (see \cite{KS}, Theorem
5.2) that  
$$
\Pm(\l;\th)=0\,, \quad \text{unless}\quad \mu\subset\l\,, \tag 2.2
$$
and that again (\cite{KS}, Corollary 4.7)
$$
\Pm(x;\th)=P_\mu(x;\th)\, + \quad\text{\sl lower degree terms} \,. \tag2.3
$$
Their proof was based on difference equations for the polynomial
$\Pm(x;\th)$. For $\th=1$ these properties follow immediately
from explicit formulas for $\Pm(x;1)$.

Explicit formulas for the polynomials $\Pm(x;\th)$, which
generalize explicit formulas for shifted Schur functions,
were found in \cite{Ok3} (see, for example, \tht{2.4} below).
In particular, they 
provide an effective proof of the existence of $\Pm(x;\th)$ and a 
different proof of \tht{2.2} and \tht{2.3}. 

We call these polynomials
{\it shifted Jack polynomials}. They are a degenerate case of shifted
Macdonald polynomials \cite{Kn,S2,Ok3}.

It easily follows from the definition that $\Pm(x;\th)$ are stable,
that is
$$
\Pm(x_1,\dots,x_n,0;\th)=\Pm(x_1,\dots,x_n;\th)
$$
provided $\mu$ has at most $n$ parts. Hence, one can consider
shifted Jack polynomials in infinitely many variables. 

The following {\it combinatorial formula} for shifted Jack polynomials
was proved in \cite{Ok3}.  Let us call a tableau
$T$ on $\mu$ a {\it reverse tableau} if its entries strictly
decrease down the columns and weakly decrease in the rows.
By $T(s)$ denote the entry in the square $s\in\mu$. We have 
$$
\Pm(x;\th)=\sum_T \psi_T(\th) \prod_{s\in\mu} 
\left(x_{\Ts}-{a'(s)}+\th\,{l'(s)}
\right) \,, \tag 2.4
$$
where the sum is over all reverse tableaux on $\mu$ with entries
in $\{1,2,\dots\}$ and $\psi_T(\th)$ is the same $\th$-weight of a tableau
that enters the combinatorial formula for ordinary Jack polynomials
(see \cite{St} or \cite{M},(VI.10.12))
$$
P_\mu(x;\th) = \sum_T \psi_T(\th) \prod_{s\in\mu} x_{\Ts} \,.  \tag 2.5
$$
The coefficients $\psi_T(\th)$ are rational functions of $\th$.

\subhead
3.~Binomial formula 
\endsubhead
Given a partition $\mu$ and a number $t$ set
$$
\align
H'(\mu)&=\prod_{s\in\mu}(a(s)+\th\, l(s)+\th)\,, \\
(t)_\mu &= \prod_{s\in\mu}(t+a'(s)-\th\, l'(s))\,.
\endalign
$$
If $\mu=(m)$ then $(t)_\mu$ is the standard shifted factorial.
We have (see \cite{M}, VI.10.20) 
$$
P_\mu(1,\dots,1;\th) = (n\th)_\mu \big/ H'(\mu) \,. \tag3.1
$$
Set (see \cite{M}, VI.10.16)
$$
Q_\mu(x;\th) = \frac{H'(\mu)}{H(\mu)}\, P_\mu(x;\th) \,.
$$
The main result of this note is the following
\proclaim{Theorem}
$$
\frac{P_\l(1+x_1,\dots,1+x_n;\th)}
{P_\l(1,\dots,1;\th)} =
\sum_\mu 
\frac{ \Pm(\l_1,\dots,\l_n;\th)\, Q_\mu(x_1,\dots,x_n;\th)}
{(n\th)_\mu}
\,. \tag3.2
$$
\endproclaim 

Note that by virtue of \tht{2.2} the summation in the binomial
formula \tht{3.2} is only over such $\mu$ that
$$
\mu\subset\l\,.
$$

Let $\O1$ be the local ring of symmetric rational functions regular
at the point
$$
\bar1=(1,\dots,1) \tag3.3
$$
and let $\m1\subset\O1$ be its maximal ideal. We call functionals
$$
\psi \in \left(\O1/\m1^k\right)^*\subset \O1^*
$$
{\it symmetric distributions} of order $\le k$ 
supported at the point $\bar1$ defined in \tht{3.3}.
Denote by $\S1$ the space of all symmetric 
distributions supported at $\bar1$. Looking
at the highest derivatives one easily proves that the map
$$
\Dn\to\S1 \tag3.4
$$
which takes an operator $D$ to the distribution $\psi_{\sssize D}$
$$
\psi_{\sssize D}(f)=(D\,f)(1,\dots,1)
$$
is an order preserving isomorphism of linear spaces. 

\demo{Proof of the theorem}
Consider the Taylor expansion of $P_\l(x;\th)$ about the point $\bar1$.
Since Jack polynomials form a linear basis in symmetric polynomials
this expansion has the form
$$
P_\l(1+x_1,\dots,1+x_n;\th) = \sum_\mu \psi_\mu(P_\l)\, P_\mu(x;\th)\,,
\quad \psi_\mu\in\S1\,.\tag3.5
$$
Note that $\psi_\mu$ is of order $\le|\mu|$. By virtue of the isomorphisms
\tht{2.1} and \tht{3.4} we have
$$
P_\l(1+x_1,\dots,1+x_n;\th) = P_\l(1,\dots,1)
\sum_\mu f_\mu(\l)\, P_\mu(x;\th)\,, \quad
f_\mu\in\Ln\,, \tag3.6
$$
where $\operatorname{deg} f_\mu\le|\mu|$. On the other hand
$$
P_\l(1+x_1,\dots,1+x_n;\th) =
P_\l(x_1,\dots,x_n;\th)\,+\,\text{\sl lower degree terms}\,. \tag3.7
$$
{F}rom \tht{3.6} and \tht{3.7} we obtain
$$
f_\mu(\l)=\cases
0,  & |\lambda|\le|\mu|,\lambda\ne\mu\,,\\
1/P_\mu(1,\dots,1), & \lambda=\mu\,. 
\endcases
$$

Since $f_\mu$ is a polynomial of degree $\le|\mu|$ it is 
completely determined by its values at the points $\l$,
$|\lambda|\le|\mu|$. Therefore,  $f_\mu$ is proportional
to  $\Pm$ and by  \tht{3.1}
$$
f_\mu(x) = \frac{H'(\mu)}{(n\th)_\mu H(\mu)} \Pm(x;\th) \,,
$$
which concludes the proof. \qed
\enddemo

\remark{\bf Remark}
The coefficients ${\binom\lambda\mu}_\th$ in the expansion
$$
\frac{P_\l(1+x_1,\dots,1+x_n;\th)}
{P_\l(1,\dots,1;\th)} =
\sum_\mu {\binom\lambda\mu}_\th
\frac{P_\mu(x_1,\dots,x_n;\th)}
{P_\mu(1,\dots,1;\th)}\,.
$$
are
called the {\it generalized binomial coefficients}. 
They were studied by
C.~Bingham \cite{Bi} in the special case $\th=1/2$, using
group--theoretical methods, and by M.~Lassalle \cite{La} for general
$\th$. Then a proof of the main results announced in \cite{La} was proposed by
J.~Kaneko (see \cite{K}, section 5).
Our theorem says that this binomial coefficient
$$
{\binom\lambda\mu}_\th=\frac{\Pm(\lambda;\th)}{H(\mu)} \tag 3.8
$$
is a shifted symmetric
polynomial in $\l$ which by \tht{2.4}
and \tht{2.5} is just as complex as the Jack
polynomial $P_\mu(x;\th)$ itself. The main results of Lassalle
can be easily deduced from \tht{3.8}. For example,
it is clear that \tht{3.8} does not
depend on $n$ and vanishes unless $\mu\subset\l$.
The recurrence relation (see \cite{La},
th\'eor\`emes 2, 3, 5 and corollaire in
section 6) is equivalent to the formula \tht{5.2} below. 

Note also that the coefficients \tht{3.8} admit a simple 
formula when $\th=1$ (see A.~Lascoux \cite{Lasc} and Example 10 in
\cite{M}, section I.3) or when $\lambda_1=\dots=\lambda_n$ (see 
J.~Faraut and A.~Kor\'anyi \cite{FK},
Prop.\ XII.1.3 (ii)). 

Binomial formulas for characters of classical groups are
discussed in \cite{OO2}
\footnote{
Recently the first author found a binomial formula for
Macdonald polynomials (both shifted and ordinary) and
also for Koornwinder polynomials, see \cite{Ok4-5}.
}.
\endremark

\subhead
4.~Bessel functions (\cite{D,Op,J})
\endsubhead
For a real vector $l=(l_1,\dots,l_n)$ denote by
$$
[l]=([l_1],\dots,[l_n])
$$
its integral part. 
Suppose that $l_1\ge\dots\ge l_n$.
By definition, put
$$
F(l,x;\th) =\lim_{\kappa\to\infty}
\frac{P_{[\kappa l]} (1+x_1/\kappa,\dots,1+x_n/\kappa;\th)}
{P_{[\kappa l]} (1,\dots,1;\th)} \,. \tag4.1
$$
{F}rom \tht{3.2} and \tht{2.3} we have
\proclaim{Proposition}
$$
F(l,x;\th) =\sum_\mu 
\frac{ P_\mu(l_1,\dots,l_n;\th)\, Q_\mu(x_1,\dots,x_n;\th)}
{(n\th)_\mu}
\,. \tag4.2
$$
\endproclaim

\demo{Proof}
We have (one easily checks absolute convergence of all
series)
$$
\align
F(l,x;\th)&=
\lim_{\kappa\to\infty}
\frac{P_{[\kappa l]} (1+x_1/\kappa,\dots,1+x_n/\kappa;\th)}
{P_{[\kappa l]} (1,\dots,1;\th)}\\
&=\sum_\mu \lim_{\kappa\to\infty}
\frac{ \Pm([\kappa l_1],\dots,[\kappa l_n];\th)\,
 Q_\mu(x_1/\kappa,\dots,x_n/\kappa;\th)}
{(n\th)_\mu} \\
&=\sum_\mu \lim_{\kappa\to\infty}
\frac{ \kappa^{-|\mu|} \Pm([\kappa l_1],\dots,[\kappa l_n];\th)\,
 Q_\mu(x_1,\dots,x_n;\th)} 
{(n\th)_\mu} \\
&=\sum_\mu 
\frac{ P_\mu(l_1,\dots,l_n;\th)\,
 Q_\mu(x_1,\dots,x_n;\th)} 
{(n\th)_\mu} \,,
\endalign
$$
where the second equality is based on  the binomial formula
and the last equality follows from \tht{2.3} \qed.
\enddemo

We call $F(l,x;\th)$ the {\it Bessel functions}. They are 
in the same relation to Jack polynomials as ordinary Bessel
functions to Jacobi polynomials. They are eigenfunctions of
the corresponding degeneration of Sekiguchi operators.
The formula \tht{4.2} makes obvious the following 
known symmetry 
$$
F(l,x;\th)=F(x,l;\th)\,.
$$

Let $H(n,\Bbb R)$ denote the space of real symmetric matrices of order $n$,
let $X,Y\in H(n,\Bbb R)$, and let $x=(x_1,\dots,x_n)$ and
$y=(y_1,\dots,y_n)$ be the spectra of $X$ and $Y$. The compact orthogonal
group $O(n)$ acts on the space $H(n,\Bbb R)$ by conjugations. One can prove
that
$$
F(y,x;1/2)=\int_{O(n)} \operatorname{exp}
(\operatorname{tr}(YuXu^{-1}))du,\tag4.3
$$
where $du$ is the normalized Haar measure on $O(n)$. (Idea of proof: the
normalized polynomials
$$
P_\lambda(z_1,\dots,z_n;1/2)/P_\lambda(1,\dots,1;1/2)
$$
are the spherical functions on the compact symmetric space $U(n)/O(n)$
while the integrals \thetag{4.3} are essentially the spherical functions
on the associated Euclidean type symmetric space 
$O(n)\ltimes H(n,\Bbb R)$. It is well--known (see \cite{DR}) that
the spherical functions of Euclidean type can be obtained from the
spherical functions of compact type by the limit transition \thetag{4.1}.)

A similar interpretation of the functions $F(y,x;\theta)$ also exists for
$\theta=1$ and $\theta=2$. Then one has to consider orbits in the spaces
$H(n,\Bbb C)$ and $H(n,\Bbb H)$ of complex and quaternionic Hermitian
matrices, respectively. An application
of the expansion \thetag{4.2} for $\theta=1$ is given in \cite{OV}; similar
results also hold for $\theta=1/2$ and $\theta=2$.

\subhead
5.~Formula for $\th$-dimension of a skew diagram $\l/\mu$
\endsubhead
Define the $\th$-dimension $\td \l/\mu$ of a skew diagram $\l/\mu$ as
the following coefficient
$$
(\sum x_i)^k P_\mu(x;\th) = \sum_{|\l|=|\mu|+k} \td\l/\mu \,\, P_\l(x;\th)\,,
\quad k=1,2,\dots \,. \tag5.1
$$
If $\th=1$ then $\td\l/\mu$ equals the number of the standard tableaux
on $\l/\mu$; for general $\th$ each tableau $T$ is counted with a certain 
weight $\psi'_T(\th)$ given in \cite{M}, section VI.6. In particular,
if $\mu=\emptyset$ then (see \tht{5.6} below for a proof using shifted Jack
polynomials)
$$
\td \l = |\l|! / H(\l) \,.
$$
We have
\proclaim{Proposition}
$$
\frac{\td \l/\mu}{\td \l} = 
\frac{\Pm(\l;\th)}{|\l|(|\l|-1)\cdots(|\l|-|\mu|+1)} \,. \tag5.2
$$
\endproclaim
This proposition can be deduced from \tht{3.8} and \cite{La,K}
(see the remark above), but it easier to give a direct proof, 
which uses the very same argument
as in \cite{OO}, section 9. 

First one proves that
\proclaim{Lemma}
$$
(\sum x_i-|\mu|) \Pm(x;\th) = \sum_{|\l|=|\mu|+1} 
\td\l/\mu \,\, P^*_\l(x;\th) \,. \tag 5.3
$$
\endproclaim
\demo{Proof} Let $f$ be the difference of the LHS and RHS in \tht{5.3}.
By \tht{2.3} and \tht{5.1} we have $\deg f\le|\mu|$. One easily checks that
$f(\nu)=0$ for all $|\nu|\le|\mu|$. It follows that $f=0$. \qed
\enddemo 
Applying this lemma $k$ times, where $k=|\l|-|\mu|$, we obtain
$$
(\sum x_i-|\mu|)\cdots(\sum x_i-|\mu|-k+1)
 \Pm(x;\th) = \sum_{|\nu|=|\mu|+k} 
\td\nu/\mu \,\, P^*_\nu(x;\th) \,. \tag5.4
$$
Evaluation of \tht{5.4} at $x=\l$ gives
$$
(|\l|-|\mu|)! \,\Pm(\l) = \td\l/\mu\,\, P^*_\l(\l)\,. \tag5.5
$$
In particular, for $\mu=\emptyset$ one has  
$$
|\l|! = \td\l\, P^*_\l(\l)\,. \tag5.6
$$
Dividing \tht{5.5} by \tht{5.6} we obtain \tht{5.2}.\qed 

\subhead 
6.~Integral representation of Jack polynomials
\endsubhead
Write $\nu\prec\l$ if
$$
\l_1\ge\nu_1\ge\l_2\ge\nu_2\ge\dots\ge\nu_{n-1}\ge\l_n \,.
$$
For any number $r$ set $(t)_r=\G(t+r)/\G(t)$. 
We have (see \cite{M}, VI.7.13${}'$)
$$
P_\l(x_1,\dots,x_n;\th)=\sum_{\nu\prec\l} \psi_{\l/\nu} 
P_\nu(x_1,\dots,x_{n-1}) x_n^{|\l/\nu|}\,, \tag6.1
$$
where
$$
\psi_{\l/\nu}=
\prod_{i\le j\le {n-1}}
\frac{
(\nu_i-\l_{j+1}+\th(j-i)+1\ft 
(\l_i-\nu_j+\th(j-i)+1\ft}
{(\l_i-\l_{j+1}+\th(j-i)+1\ft
(\nu_i-\nu_j+\th(j-i)+1\ft}
\,. \tag6.2
$$
Suppose $\mu$ has less than $n$ parts.
The formulas \tht{3.2} and \tht{6.1} together imply 
$$
\Pm(\l) = \frac{(n\th)_\mu}{((n-1)\th)_\mu} \sum_{\nu\prec\l}
\psi_{\l/\nu} 
\frac{P_\nu(1,\dots,1)}{P_\l(1,\dots,1)}\, \Pm(\nu) \tag6.3
$$
where $1$ is repeated $n$ times in the denominator and 
$(n-1)$ times in the numerator. Note that after substitution of \tht{6.2}
and of \tht{3.1} rewritten as 
$$
P_\l(1,\dots,1)=\prod_{i<j\le n} (\l_i-\l_j+\th(j-i))_\th \prod_{i\le n}
\G(\th)/\G(\th i) \tag6.4
$$
the formula \tht{6.3} becomes the formula \tht{7.16} from \cite{Ok3}. 
We replace $\l$ and $\mu$ in \tht{6.3} by their multiples $\kappa \l$ and
$\kappa \mu$ and let $\kappa\to\infty$. By \tht{2.3}
$$
\Pm(\kappa\l)/\kappa^{|\mu|} \to P_\mu(\l)\,, \quad \kappa\to\infty \,.
$$
Introduce the following product of beta-functions
$$
C(\mu,n)=\prod_{i\le n-1} B(\mu_i+(n-i)\th,\th)\,.
$$
Set
$$
\Pi(\l,\nu;\th)=\prod_{i\le j}(\l_i-\nu_j)^{\th-1}
\prod_{i>j}(\nu_j-\l_i)^{\th-1}\,.
$$
Using the well known relation $(t)_\th/t^\th\to 1$, $t\to\infty$,
one obtains from \tht{6.3} the following 
\proclaim{Proposition}  
If $\mu$ has $<n$ parts then
$$
P_\mu(\l)=
\frac 1{C(\mu,n)} \frac1{V(\l)^{2\th-1}}
\int_{\l_2}^{\l_1}d\nu_1\!\dots\!
\int_{\l_n}^{\l_{n-1}}\!d\nu_{n-1}\,
 P_\mu(\nu)\,V(\nu)
\, \Pi(\l,\nu;\th)\,.
$$
\endproclaim
Here $\l_1\ge\dots\ge\l_n$ are arbitrary real (this assumption
is not essential since $P_\mu$ is symmetric). 
If $\th=1,2,\dots$ then
the integrand is holomorphic and $\l$ can be arbitrary complex.
Iteration of this proposition together with \tht{1.2} gives an
integral representation for all Jack polynomials. 

This proposition was first obtained by G.~Olshanski (unpublished) for
special values $\th=1/2,1,2$ by using the group theoretic
interpretation of the corresponding Jack polynomials. Then
A.~Okounkov \cite{Ok3} found a general method which works for any $\th$ and
even for Macdonald polynomials. The proof presented above differs
from that from \cite{Ok3}.  

\subhead
7.~Other applications
\endsubhead
The binomial theorem itself and the formula \tht{5.2} result in
an $\th$-analog of the Vershik-Kerov theorems \cite{VK1,VK2} about
characters of $U(\infty)$ and $S(\infty)$. 
These results will be discussed elsewhere.

\end